\documentclass{ifacconf}

\usepackage{graphicx}      % include this line if your document contains figures
\usepackage{natbib}        % required for bibliography
\usepackage{amsmath}
\usepackage{amsfonts}

\usepackage{setspace}
\usepackage{booktabs}
\usepackage{multirow}
\usepackage{algorithm}
\usepackage[noend]{algpseudocode}
\begin{document}
\begin{frontmatter}

%\title{A Grid-Based Bayesian State Observer with Markov Chain Monte Carlo Propagation Stage \thanksref{footnoteinfo}} 

\title{An approximation of the Bayesian state observer with Markov chain Monte Carlo propagation stage \thanksref{footnoteinfo}} 

\thanks[footnoteinfo]{This work has been supported by the COMET-K2 Center of the Linz Center of Mechatronics (LCM) funded by the Austrian federal government and the federal state of Upper Austria.
}

\author[First]{L. Ecker} 
\author[First]{K. Schlacher} 
%\author[Third]{Third C. Author}

%Johannes Kepler University Linz, Institute of Automatic Control and Control Systems Technology, Altenberger Str. 69, A-4040 Linz, Austria, e-mail: lukas.ecker@jku.at
\address[First]{Institute of Automatic Control and Control Systems Technology, Johannes Kepler University Linz,  Altenberger Str. 69, A-4040 Linz, Austria (e-mail: \{lukas.ecker, kurt.schlacher\}@jku.at).}
%\address[Second]{Johannes Kepler University Linz, Institute of Automatic Control and Control Systems Technology, Altenberger Str. 69, A-4040 Linz, Austria(e-mail: author@lamar. colostate.edu)}
%\address[Third]{Electrical Engineering Department,Seoul National University, Seoul, Korea, (e-mail: author@snu.ac.kr)}

\begin{abstract}                
The state estimation problem for nonlinear systems with stochastic uncertainties can be formulated in the Bayesian framework, where the objective is to replace the state completely by its probability density function.
Without the restriction to selected system classes and disturbance properties, the Bayesian estimator is particularly interesting for highly nonlinear systems with non-Gaussian noise.
%Since there are hardly any analytical solutions of the general Bayesian observer, various approximation methods have been developed for the numerical implementation of the state estimation filter.
The main limitations of Bayesian filters are the significant computational costs and the implementation problems for higher dimensional systems.  
The present paper introduces a piecewise linear approximation of the Bayesian state observer with Markov chain Monte Carlo propagation stage and kernel density estimation. 
These methods are suitable for the prediction of multivariate probability density functions.
The piecewise linear approximation and the proposed algorithms can increase the estimation performance at reasonable computational cost.
The estimation performance is demonstrated in a benchmark comparing the Bayesian state observer with an extended Kalman filter and a particle filter.

%Die grundlegendste Beschreibung des zustandsschätzsproblem für nichtlineare systeme mit nicht gausschem Rauschen ist gegeben durch die Theorie von Bayes.
	
%These instructions give you guidelines for preparing papers for the $10^{\mathrm{th}}$ Vienna Conference on Mathematical Modelling. Please use this document as a template to prepare your manuscript. For submission guidelines, follow instructions on paper submission system as well as the event website. 

%Full contributions for the $10^{\mathrm{th}}$ Vienna Conference on Mathematical Modelling are limited to 6 pages. Please keep the abstract of your paper within a limit of approximately 300 words.
\end{abstract}


\end{frontmatter}
%===============================================================================

\section{Introduction}

Many processes from a wide range of scientific fields are modelled as nonlinear stochastic dynamic systems.
Among other things, a mathematical model offers the possibility of a system-theoretical investigation. 
In this way, for example, influences can be estimated, predictions can be made, and trajectories or control strategies for system operation can be designed.
For these purposes, knowledge of the entire system state is often necessary.
In some cases, however, it is not possible to measure all of the state variables associated with the model, whether for technical, environmental, or economic reasons.
In addition, the available measurements may be superimposed by non-deterministic influences such as measurement noise or process disturbances. For the design and implementation of a state observer, sampled-data systems of continuous-time processes or inherent discrete-time systems are considered.
%Therefore, various filters are used to extract information from partly uncertain environments in order to obtain estimates of non-measurable or poorly measurable quantities.
%These state observers differ in their theoretical approach, complexity, implementation effort, and computational intensity depending on the class of systems under consideration.
% Abstand am Besten
A fundamental formulation of the nonlinear state estimation problem is given in the Bayesian framework.
The idea of the Bayesian state observer (BSO) is to describe the estimation of the state completely by its probability density function (pdf).
The algorithm consists of two stages, one propagating the density to the next time step and the other inferring additional state information through measurements. 
In its most general form, i.\,e., for the nonlinear and non-Gaussian noise scenario, the prediction step corresponds to the Chapman–Kolmogorov equation and the update step complies with Bayesian inference.
An analytically computable solution in closed-form is generally hardly available.
Therefore, several approaches are used to approximate the Bayesian state estimation problem towards a numerical implementation.
The increased computational effort compared to the extended Kalman filter (EKF), for example, is a trade off for improved estimation performance.
% Abstand am Besten
In this paper, a grid-based approximation of the BSO with multilinear interpolation of the pdfs on a finite set of grid nodes is presented. While in this scenario the inference step can be efficiently performed by a multi-threaded algorithm, the propagation step is accomplished by a Markov chain Monte Carlo (MCMC) method, where a large number of samples distributed according to the a posteriori density of the state are propagated through the dynamic equation and reconstructed into an a priori pdf using kernel density estimation (KDE).

Therefore, the paper is structured as follows. In the remaining part of the introduction a brief survey of the Bayesian filtering theory is given, where different known numerical approximations are referenced. In addition, the notation used from estimation and probability theory is briefly explained. In Section \ref{sec:bayesian_state_observer}, the fundamental formulation of the Bayesian state observer is introduced and the two main stages of the filter are explained. In Section \ref{sec:numerical_approximation}, a grid-based observer with Markov chain Monte Carlo propagation stage and kernel density estimation is presented as numerical approximation of the Bayesian state estimation problem. Finally, experimental results of the proposed Bayesian approximation for a benchmark example are presented in Section \ref{sec:experimental_results}.

\subsection{Bayesian filtering theory}
A general approach to the state estimation problem for nonlinear stochastic dynamic systems was formulated more than four decades ago in the Bayesian framework, see, e.\,g., \cite{HoLee:64}.
The analytical equations of the Bayesian state observer, described e.\,g., in \cite{Ristic:2003}, form a basic unified framework for many of the well-known filters, such as the Kalman filter and the particle filter, see \cite{Sim:06}. 
%While the Kalman filter can be derived considering a linear system and Gaussian noise, the particle filter was invented to implement the Bayesian estimator numerically, see \cite{Sim:06}.
%As mentioned, in the Bayesian approach, a generally infinite-dimensional probability density function is used for a complete statistical representation of the state, see \cite{Papoulis:2002}. 
However, in the most general case of the Bayesian state estimation problem, there is often no exact parametric solution. %Therefore, various approximations of the BSO exist that allow for numerical implementation.
%The objective of the Bayesian approach, is to recursively construct the complete probability density function, which is generally infinite-dimensional, and to use it for the statistical representation of the state, see \cite{Papoulis:2002}.
% However, the analytical equations for the processing of the density contain multi-dimensional integration terms for which the existence of a closed solution is not certain and the computational effort would also be hardly manageable, see \cite{KRAMER:1988}. 
%T most general case of the Bayesian state estimation problem often yields no exact parametric solution.
%Therefore, when considering the general non-linear state estimation problem for which no exact solution can be found, numerical approximations are required. 
An extensive overview of Bayesian filtering theory can be found, e.\,g., in \cite{Chen:03}, where the author further discusses a wide range of different numerical approximations.
The BSO approximations can be roughly divided into grid-based methods, point-mass methods, parametric methods and non-parametric methods, see \cite{Stano:2013}.
While grid-based and point-mass methods typically involve approximation of the state space, see, e.\,g., \cite{KRAMER:1988}, parametric and non-parametric methods look for a solution to the estimation problem in a finite-dimensional parameter space or are based on the idea of Monte Carlo integration.
Among the best-known non-parametric methods are the various variants of the particle filter, which also include bootstrap filters, see \cite{Gordon:93}, and sequential Monte Carlo filters, see \cite{Andrieu:2001}.
%These nonlinear estimators restrict the statistical representation to a set of random particles, and often do not consider the complete pdf of the state.

 %such as Particle filters, see also \cite{Sim:06}, which also include bootstrap filters, see \cite{Gordon:93}, and sequential Monte Carlo filters, see \cite{Andrieu:2001}, are among the best-known non-parametric methods. 

 %These non-linear estimators restrict the statistical representation to a set of random state vectors, i.e. the particles, and do not handle the complete pdf of the state. 

\subsection{Mathematical notation}

The notation and main formulas from probability and estimation theory are briefly mentioned. A detailed description can be found, e.\,g., in \cite{Papoulis:2002} and \cite{Sim:06}.
%Since processes with non-deterministic influences are considered in the following, the concept of random variables and their probability is important. Nevertheless, this introduction is limited only to the basics, since a complete and rigorous definition falls within the scope of measure and integration theory, see \cite{Luenberger:1997}.
The system state and the disturbances are interpreted as tuples of continuous real-valued random variables. Given a  multivariate random variable $X$ with domain $\mathcal{D}$, the monotonically increasing probability distribution function, i.\,e., the probability that $X$ is less than or equal to a vector $x = [x_1,\ldots,x_n]^T \in \mathcal{D}$, is denoted by $P_X(x) = \Pr(X\le x) \in \left[0,1\right]$. Its corresponding multivariate probability density function is defined as the derivative $p_X(x) = \frac{\partial}{\partial x} P_X(x)= \frac{\partial^n}{\partial x_1,\ldots,\partial x_n} P_X(x_1,\ldots,x_n) \ge 0$ with the property  $\int_{\mathcal{D}} p_X(x) \mathrm{d}x = 1$.
The subscripts are omitted when the associated random variables are evident from the context, i.\,e., the pdfs of the two multivariate random variables $X$ and $V$ are denoted $p(x)$ and $p(v)$, respectively.
Moreover, the conditional pdf of $X$ given $V$ is interpreted as
$p(x|v) = \frac{\partial}{\partial x} P(x|v) = \frac{\partial}{\partial x} \Pr(X \le x|V = v)$.
Given a diffeomorphic mapping $g: \mathcal{D}_X \rightarrow \mathcal{D}_V, x \mapsto g(x) = v $, the statistics of the random variable $V = g(X)$ can be expressed in terms of the statistics of X, i.\,e.,
\begin{equation}\label{equ:pdf_transformation}
p(v) = \left(p(x) \det(\partial_x g)^{-1}\right) \circ g^{-1}(v)
\end{equation}
where $\partial_x g$  denotes the Jacobian matrix of the transformation $g$, see, e.\,g., \cite{Papoulis:2002}.
A connection between the joint statistics of $X$ and $V$, i.\,e., $p(x,v)$, and its conditional pdf $p(x|v)$ and $p(v|x)$, respectively, is given by the total probability $p(v) = \int_{\mathcal{D}_X} p(v|x)p(x) \mathrm{d}x$ and Bayes' theorem according to
\begin{equation} \vspace*{-0.0cm}
p(x|v) = \frac{p(x,v)}{p(v)} = \frac{p(v|x)p(x)}{\int_{\mathcal{D}_X} p(v|x)p(x) \mathrm{d}x}.
\end{equation}
The estimates or approximations of stochastic variables are noted with the hat symbol, e.\,g., $\hat{x}$. Moreover, considering sampled-data systems with $k$ denoting a particular point in time, estimates can be distinguished into a priori and posteriori quantities, e.\,g., $\hat{x}_k^-$ and $\hat{x}_k^+$.
The phrase a priori indicates that the estimate $\hat{x}_k^-$ takes all measurements $Y_{k-1} = \{y_0, \ldots, y_{k-1}\}$ up to time $k-1$ into account. 
In contrast, the a posteriori term designates an estimate $\hat{x}_k^+$ with measurements $Y_k = \{y_0,\ldots,y_k\}$ up to time $k$.

%In estimation the
%Next we see a few subsections.
%
%
%\subsection{Terms from probability theory}
%\subsection{Distinction in a priori and a posteriori}
%
%
%\subsection{Markov Processing}
%
%
%\begin{enumerate}
%	\item pdf
%	\item transformation pdf
%	\item Bayes Theorem
%	\item A priori a posteriori
%	
%	\item  Markov chain Monte Carlo method
%\end{enumerate}
%
%\subsection{Fundamental rules}
%

\section{Bayesian State Observer}\label{sec:bayesian_state_observer}

The intention is to design a state observer for the nonlinear system described by the state and measurement equations
\hspace*{-0.5cm}
\begin{subequations}\label{equ:nonlinear_model}
	\begin{align}
		x_{k+1} = f_k(x_k,u_k,w_k) \label{equ:nonlinear_system_state_equation}\\
		y_{k} = h_k(x_k,u_k,v_k)
		\label{equ:nonlinear_system_output_equation}
	\end{align}
\end{subequations}
with index $k$ denoting the current time, the system state $x_k \in \mathbb{R}^{n_x}$,  the system input $u_k \in \mathbb{R}^{n_u}$ and output $y_k \in \mathbb{R}^{n_y}$, the random process disturbance or noise $w_k \in \mathbb{R}^{n_x}$ and the random measurement noise $v_k \in \mathbb{R}^{n_y}$.
The model (\ref{equ:nonlinear_model}) corresponds to a stochastic process with random variables $x_0$, $w_k$, $v_k$ and deterministic input $u_k$, which implies that the state $x_k$ is also a random variable.
The actual values of the input $u_k$ and the output $y_k$ at the current time $k$ are known, while those of the non-deterministic noises $w_k$ and $v_k$ are not available. 
Instead, only knowledge of the probability density functions $p(w_k)$ and $p(v_k)$ is assumed.
The same applies for the initial value of the state $x_k$, which is represented by its density $p(x_0)$.
After the initialization process, the Bayesian state observer algorithm is divided into two steps that are executed alternately. As with the Kalman filter, these steps are referred to as prediction and update.
%The algorithms and implementations are described in detail in the Subsections (\ref{subsec:bayesian_prediction}) and (\ref{subsec:bayesian_inference}), although the associated basic ideas will be briefly outlined here. 
The prediction or propagation step intends to estimate the a priori state $\hat{x}^-_{k}$ at the time $k$ given the best estimation of the a posteriori state $\hat{x}^+_{k-1}$ of the previous time step $k-1$.
On the other hand, the update or inference step incorporates a received measurement $y_k$ into its estimate, i.\,e., it takes the a priori state $\hat{x}^-_{k}$ and infer the current measurement $y_k$ to obtain the a posteriori state $\hat{x}^+_k$. The Bayesian state observer does not restrict its estimations to a few stochastic moments. Rather, the descriptive pdf of the state is propagated from one time step to the next, also inferring measurements obtained.
At each time step, various methods can be used to approximate an optimal estimate $\hat{x}^+_k$ from an a posteriori pdf. 
%As an initial condition for the Bayesian state observer, knowledge of the initial state probability density $p_(x_0)$ and the probability densities of the
%perturbations $p(w_k)$ and $p(v_k)$ for $k = 0,1,\ldots $ is required.
\subsection{Bayesian Prediction}
\label{subsec:bayesian_prediction}
Given the probability density of the state at time $k-1$, i.\,e., $p(x_{k-1}|Y_{k-1})$, the objective of the prediction or propagation step within the Bayesian state observer algorithm is to predict the conditional pdf of $x_{k}$ given all the measurements prior to time $k$, i.\,e., $p(x_{k}|Y_{k-1})$. 
An a priori estimate $\hat{x}^-_{k}$ can be abstracted from the conditional pdf $p(x_{k}|Y_{k-1})$ according to various optimality criteria. For instance, the estimate of the maximum a posteriori probability (MAP) is given by
\begin{equation}
\hat{x}^-_k = \underset{x_k \in \mathcal{D}}{\arg \max }\; p(x_k|Y_{k-1}).
\end{equation}
Nevertheless, the expectation value as in the Kalman filter would also be conceivable, although caution is required with bimodal probability densities.
However, even the MAP estimate of a bimodal distribution in which the highest mode is uncharacteristic of most of the distribution can have dramatic effects on subsequent estimates.
The analytical propagation equation is obtained as marginal probability density function of $p(x_k,x_{k-1}|Y_{k-1})$, see (\ref{equ:prediction_marginal_density_1}), where Bayes' Rule allows to rewrite the density as $p(x_k,x_{k-1}|Y_{k-1}) = p(x_k|x_{k-1},Y_{k-1})p(x_{k-1}|Y_{k-1})$. The Markov process characteristics of (\ref{equ:nonlinear_system_state_equation}), i.\,e., the probability of the current state given the immediately previous one is conditionally independent of earlier states as well as outputs, yields $p(x_k|x_{k-1},Y_{k-1}) = p(x_k|x_{k-1})$.
Therefore, a recursive assignment is given by 
\begin{align}
p(x_k|Y_{k-1}) &= \int_\mathcal{D} p(x_k,x_{k-1}|Y_{k-1}) \mathrm{d}x_{k-1} \label{equ:prediction_marginal_density_1}\\
			   &= \int_\mathcal{D} p(x_k | x_{k-1}) p(x_{k-1}|Y_{k-1}) \mathrm{d}x_{k-1},\label{equ:prediction_marginal_density_2}
\end{align}
where $p(x_{k-1}|Y_{k-1})$ is the known previous a posteriori pdf and $p(x_k|x_{k-1})$ is the pdf of $x_k$ given a state $x_{k-1}$. This density is known, since it is derived from the given probability density $p(w_k)$ and the state equation (\ref{equ:nonlinear_system_state_equation}). For this purpose, we assume that the mapping (\ref{equ:nonlinear_system_state_equation}) with fixed state $x_k$ and input $u_k$ is diffeomorphic and denote the associated inverting map by $w_{k+1} = f_k^*(x_k,u_k,x_{k+1})$.
The pdf $p(x_k|x_{k-1})$, see also (\ref{equ:pdf_transformation}), is determined by 
\begin{equation*}
	p(x_k|x_{k-1}) = \left(p(w_k)\det(\partial_w f_{k-1})^{-1}\right)  \circ f_{k-1}^*(x_{k-1},u_{k-1},x_{k}),
\end{equation*}
where $\partial_{w} f_{k-1}$ corresponds to the Jacobian matrix with entries $[\partial_{w} f_{k-1}]_{ij}= \frac{\partial f_{i,k-1}}{\partial w_{j,k-1}}$.
\subsection{Bayesian Inference}
\label{subsec:bayesian_inference}
Given the current estimate at time $k$, i.e., the conditional pdf $p(x_k|Y_{k-1})$, the objective of the inference step is to infer the received output $y_k$ to obtain the a posteriori pdf $p(x_k|Y_k)$.  
%As discussed previously a estimation of the state can be given by the MAP as $$
%	\hat{x}^+_k = \underset{x_k \in \mathcal{D}}{\arg \max }\; p(x_k|Y_{k}).$$	
The MAP state estimation is determined as $$
\hat{x}^+_k = \underset{x_k \in \mathcal{D}}{\arg \max }\; p(x_k|Y_{k}).$$
	%According to Bayes' theorem the posteriori pdf can be deduced as
From the definition of conditional density, see, e.\,g., \cite{Papoulis:2002}, the posteriori pdf
\begin{equation*}
	p(x_k|Y_k) = \frac{p(x_k,Y_k)}{p(Y_k)}= \frac{p(x_k,y_k,Y_{k-1})}{p(y_k,Y_{k-1})}
\end{equation*}
can be deduced according to Bayes' theorem and after simplification as
\begin{equation*}
%\label{equ:bayesian_inference}
p(x_k|Y_k) = \frac{p(y_k|x_k,Y_{k-1})p(x_k|Y_{k-1})}{p(y_k|Y_{k-1})},
\end{equation*} 
which is with the Markov property of (\ref{equ:nonlinear_model}) equivalent to
\begin{equation}
\label{equ:bayesian_inference}
p(x_k|Y_k) = \frac{p(y_k|x_k)}{p(y_k|Y_{k-1})}p(x_k|Y_{k-1}).
\end{equation} 
	The Bayesian inference (\ref{equ:bayesian_inference}) provides a relationship between the a priori $p(x_k|Y_{k-1})$ and the a posterior density $p(x_k|Y_k)$. While the denominator of (\ref{equ:bayesian_inference}) corresponds to the normalization factor $$p(y_k|Y_{k-1}) = \int_{\mathcal{D}} p(y_k|x_k)p(x_k|Y_{k-1}) \mathrm{d}x_k,$$ the numerator $p(y_k|x_k)$ weights the likelihood of the current measurement $y_k$. The numerator can be determined by the output equation (\ref{equ:nonlinear_system_output_equation}) and the pdf of the measurement noise $p(v_k)$. The diffeomorphic inverting map of (\ref{equ:nonlinear_system_output_equation}) with fixed state $x_k$ and input $u_k$ is denoted by $v_{k} = h^*_k(x_k,u_k,y_k)$.
	The pdf $p(y_k|x_k)$ is thus given by 
%	\begin{equation}
%	p(y_k|x_k) = \left(p(v_k) \det(\mathrm{J_g})\right) \circ g_k(x_k,u_k,y_k),
%	\end{equation}
%	where $\mathrm{J_g}$ denotes the Jacobian matrix, whose $(i,j)$-th entry is given by $[\mathrm{J_g}]_{ij}= \frac{\partial g_{i,k}}{\partial y_{j,k}}$.
	\begin{equation}
	\label{equ:bayesian_inferance_measurement_liklihood}
	p(y_k|x_k) = \left(p(v_k) \det(\partial_{y} h_k)^{-1}\right) \circ h^*_k(x_k,u_k,y_k),
	\end{equation}
	where $\partial_{y} h_k$ denotes the Jacobian matrix, whose $(i,j)$-th entry is given by $[\partial_{y} h_k]_{ij}= \frac{\partial h_{i,k}}{\partial y_{j,k}}$.

	%\begin{equation}
	%\mathrm{J_g} = \frac{\partial(g_{1,k},\ldots,g_{n_y,k})}{\partial (y_{1,k},\ldots,y_{n_y,k})}
	%\end{equation}

\section{Approximation and Implementation}\label{sec:numerical_approximation}
%Algorithm And Implementation

An exact solution to the Bayesian observation problem is available only for a few exceptions. Therefore, several approximation approaches have been proposed to obtain capable algorithms with improved estimation performance and acceptable computational cost, see, e.\,g., \cite{Chen:03}.
In particular, the multidimensional integration terms of the propagation stage, see (\ref{equ:prediction_marginal_density_2}), are a thorn in the flesh for the implementation of the state observer.
In this section, a grid-based piecewise linear approximation of the Bayesian state observer is presented using a Markov chain Monte Carlo method along with a kernel density estimation algorithm that replaces the usual propagation stage. The advantage of this approach is that the algorithms can be easily extended to multidimensional systems and the increasing computational effort is kept within limits.

%\begin{figure}
%	\begin{center}
%		\includegraphics[width=8.4cm]{bayesian_propagation}    % The printed column width is 8.4 cm.
%		\caption{Bifurcation: Plot of local maxima of $x$ with damping $a$ decreasing} 
%		\label{fig:bifurcation}
%	\end{center}
%\end{figure}

\subsection{Approximation of the State Space}
\label{sec:approximation_of_state_space}
In the proposed implementation of the Bayesian state observer, the $n$-dimensional state space $\mathcal{D}$ is approximated on a finite set of grid nodes, where the number of elements $N_{G,i}$ for each dimension with $i \in \{1,\ldots,n\} = \mathcal{N}_G$ serves as hyperparameter. The simple objective of this approach is to construct a grid on $\mathcal{D}$ in which the pdf of the state can be linear interpolated piecewise. The approximated one-dimensional state space for the $i$-th element of the n-dimensional state $x_k = [x_{1,k},\ldots,x_{n,k}]^T$ in time $k$ is denoted by $S_{i,k}$ and corresponds to a regular ordered set  with $N_{G,i}$ elements of real-valued grid positions between a defined upper and lower bound. The grid points in the set $\mathcal{S}_{i,k} = \{x^1_{i,k},\ldots,x^{N_{G,i}}_{i,k}\}$ do not need to be uniformly spaced and can be rasterised in each time step. Rasterization is necessary when the main feature of the pdf drifts closer to an edge of the set, the density already extends across the set, or the resolution of the pdf needs to be increased in a certain area. However, it must be ensured that the probability densities remain normalised by the approximation and the finite grid at each time step.
Nevertheless, the problem of grid design for point-mass approaches is addressed in \cite{Simandl:2002}, where also an algorithm for determining the minimum sufficient number of grid points is developed.
The approximated state space $\mathcal{D}_k$ at time $k$ is composed of the finite sets $\mathcal{S}_{i,k}$ as
$$\prod\limits_{i = 1}^{n} \mathcal{S}_{i,k}  = \{(x_{1,k},\ldots,x_{n,k})|\,x_{i,k} \in \mathcal{S}_{i,k} \}  \subset \mathcal{D}.$$
A state in the approximated state space $\mathcal{D}_{k}$  can be interpreted as an $n$-dimensional array that is uniquely indicated by an $n$-dimensional tuple $\alpha \in \mathcal{A} = \{(\alpha_1,\ldots,\alpha_n)| \;  \alpha_i \in \mathbb{N}^+ \wedge \alpha_i \leq N_{G,i} \}$ with $$x^{\alpha}_k = (x_{1,k}^{\alpha_1},\ldots, x_{n,k}^{\alpha_n}) \in \prod\limits_{i = 1}^{n} \mathcal{S}_{i,k}.$$ 
A finite pdf approximation of, for instance, $p(x_k|Y_k)$ is symbolized with the tilde accent by $\tilde{p}(x_k|Y_k)$.
The density $\tilde{p}(x_k|Y_k)$ is completely characterized by the set of values 
\begin{equation}
\label{equ:approximated_coefficients}
m^\alpha_{k|k} = \tilde{p}(x^\alpha_k|Y_k)=  p(x^\alpha_k|Y_k)
\end{equation}
for all $\alpha \in \mathcal{A}$. On the domain $x_k \in \mathcal{D} \setminus \mathcal{D}_{k}$, the density $\tilde{p}(x_k|Y_k)$ is evaluated as the $n$-linear interpolation of the coefficients $m^\alpha_{k|k}$ with $\alpha \in \partial \mathcal{A}(x_k)$ corresponding to the immediate grid node neighbours of $x_k$, see, e.\,g., \cite{Press:2007}.  The estimated densities of the a posteriori and a priori state are approximated in this way. The interpolation of the process noise $p(w_k)$ and the measurement noise $p(v_k)$ is not necessary, as they are assumed to be analytically given. The grid of the approximated state space $\mathcal{D}_k$ consists of $N_{c} =\prod_{i=1}^{n}N_{G,i}-1$ cells, where the center of the $i$-th cell is denoted as $x^i_{c}$. Therefore, the integration of a pdf can be numerically approximated by the trapezoidal rule, which is equal to the sum of the volume $\Delta x_k^i$ of the $i$-th cell times the interpolated center density $\tilde{p}(x_k|Y_k)$  at $x_c^i$  over all $N_{c}$ cells,
 \begin{equation}
 \label{equ:integral_approximation}
 \int_{\mathcal{D}} p(x_k|Y_k) \mathrm{d} x_k \approx \sum_{i=1}^{N_{c}} \tilde{p}(x^i_c|Y_k) \Delta x^i_k.
 \end{equation}

\subsection{Approximation of the Propagation Stage}
At the beginning of the propagation state at time $k$, the a posteriori density $p(x_{k-1}|Y_{k-1})$ is given.
The purpose of this step is to propagate the density with the dynamics of the considered model to predict the a priori pdf $p(x_k|Y_{k-1})$.
A direct approximation of (\ref{equ:prediction_marginal_density_2}) by computing the probability density and numerical integration on the constructed grid, as proposed in Section \ref{sec:approximation_of_state_space}, involves an enormous computational effort.
Therefore, the idea is to represent the a posteriori density by a large number of samples that are individually  propagated forward through the state equation and from which the a priori density can be estimated. 
The proposed algorithms are called Markov chain Monte Carlo and Kernel Density Estimation.

\subsection{Markov chain Monte Carlo Propagation}

The objective is to generate two sequences of random vectors distributed according to the pdfs of the a posteriori state $\tilde{p}(x_{k-1}|Y_{k-1})$ and the process noise $p(w_ {k-1})$.
Theses samples, i.,e, $x_{k-1}^{(i)}$ and $w_{k-1}^{(i)}$, are propagated pairwise by the state equation (\ref{equ:nonlinear_system_state_equation}) to obtain a third sequence $x_k^{(i)}$ representing the a priori state of the current time step. The density of the predicted a priori state $\tilde{p}(x_k|Y_{k-1})$ is reconstructed from this sequence by kernel density estimation, see Section \ref{sec:kernel_density_estimation}.
For the purpose of generating a sampled sequence of in general $n$-dimenisonal random variables of a multivariate probability density function Markov chain Monte Carlo methods are considered. Starting from an arbitary value $x_{k-1}^{(0)}$, the elements of a Markov chain $(X_{k-1}^{(t)}) = (x_{k-1}^{(0)},x_{k-1}^{(1)},x_{k-1}^{(2)},\ldots)$ describe a stochastic process, where the probability of the element $x_{k-1}^{(t+1)}$ solely depends on the state attained in the previous step $x_{k-1}^{(t)}$, i.\,e., $\Pr(x_{k-1}^{(t+1)} = x_{k-1} | (X_{k-1}^{(t)})) = \Pr(x_{k-1}^{(t+1)} = x_{k-1} | x_{k-1}^{(t)})$. 
Therefore, a Markov process is defined by a transition kernel $Q(x_{k-1}'|x_{k-1}^{(t)})$ that proposes the next sample $x_{k-1}'$. The distribution of the Markov chain converges to a stationary density $\pi(x_{k-1})$ when the existence and the uniqueness of stationary distribution is given, see \cite{Robert:2005}. These two conditions are sufficiently satisfied if the transition is detailed balanced and the Markov process is ergodic. Algorithms that construct this type of ergodic Markov chain with stationary distribution $\pi(x_{k-1})$ are called Markov chain Monte Carlo methods. This includes the well-known Metropolis-Hastings algorithm, which is often used for the sampling of high-dimensional distributions. Detailed balanced, i.\,e., the reversibility of the transition $Q(x_{k-1}'|x_{k-1}^{(t)}) \pi(x_{k-1}^{(t)}) = Q(x_{k-1}^{(t)}|x_{k-1}') \pi(x_{k-1}')$ is ensured by seperating the transition kernel $Q(x_{k-1}'|x_{k-1}^{(t)}) = q(x_{k-1}'|x_{k-1}^{(t)}) A(x_{k-1}',x_{k-1}^{(t)})$ in a proposal distribution $q(x_{k-1}'|x_{k-1}^{(t)})$ and a acceptance distribution $A(x_{k-1}',x_{k-1}^{(t)})$ with Metropolis choice
\begin{equation}\label{equ:acceptance}
A(x_{k-1}',x_{k-1}^{(t)}) = \min \left(1, \frac{\pi(x_{k-1}')}{\pi(x_{k-1}^{(t)})}
\frac{q(x_{k-1}'|x_{k-1}^{(t)})}{q(x_{k-1}^{(t)}|x_{k-1}')}\right).
\end{equation}
The proposal distribution $q(x_{k-1}'|x_{k-1}^{(t)})$, like the transition kernel $Q(x_{k-1}'|x_{k-1}^{(t)})$, proposes a next sample $x_{k-1}' $, which, depending on the acceptance probability $A(x_{k-1}',x_{k-1}^{(t)})$ will be accepted or rejected. The proposal density can be chosen arbitrarily. Assuming that the density is symmetric and therefore the second fraction in (\ref{equ:acceptance}) is omitted, the algorithm is also simply called Metropolis.
A sketch of the proposed Metropolis-Hastings method for the sampling from $\tilde{p}(x_{k-1}|Y_{k-1})$ is given in Algorithm 1.
% https://en.wikipedia.org/wiki/Metropolis%E2%80%93Hastings_algorithm
\begin{algorithm}
	\caption{Metropolis-Hastings}
	\begin{algorithmic}[1]
		     
		\State Set initial state $x_{k-1}^{(0)}$
		\State Set index variable $i = 0$
		\While{$i < N_{s}$}
			\State{Generate random variable:}
			\State{$x_{k-1}' \leftarrow q(x_{k-1}'|x_{k-1}^{(t)})$ }
			\State{Evaluate acceptance probability:}
			\State{$A(x_{k-1}',x_{k-1}^{(i)}) = \min\left(1, \frac{\tilde{p}(x_{k-1}'|Y_{k-1})}{\tilde{p}(x_{k-1}^{(i)}|Y_{k-1})}\frac{q(x_{k-1}'|x_{k-1}^{(i)})}{q(x_{k-1}^{(i)}|x_{k-1}')}\right)$}
			\State{Generate uniform random variable:}
			\State{$u \leftarrow \mathcal{U}_{[0,1]} $}
			\State{Accept or reject:}
			\If{$u \leq A(x_{k-1}',x_{k-1}^{(i)})$}
			\State{$x_{k-1}^{(i+1)} \leftarrow x_{k-1}'$}
			\Else
			\State{$x_{k-1}^{(i+1)} \leftarrow x_{k-1}^{(i)}$}
			\EndIf
			\State{Increment index variable:}
			\State{$i \leftarrow i +1$}
		\EndWhile

	\end{algorithmic}
\end{algorithm}

\subsection{Kernel Density Estimation}\label{sec:kernel_density_estimation}
After the propagation of the random samples drawn from the a posteriori distribution, the probability density function of the propagated state can be obtained by kernel density estimation. This non-parametric method takes the previous obtained finite data set of random sample vectors for a approximation of the a posteriori density $p(x_k|Y_k)$.
The kernel density estimator of $p(x_k|Y_k)$ is given as
\begin{equation}
p(x_k|Y_k) \approx \frac{1}{n \sqrt{\det({H})}} \sum_{i = 1}^{N_{s}} K\left( {H}^{-1/2} \left(x_k - x^{(i)}\right)\right)
\end{equation}
where $H \in \mathbb{R}^{n\times n}$ corresponds to the positive definite and symmetric bandwidth matrix and $K(\cdot)$ to the multivariate kernel function. In common, the Epanechnikov and normal functions are used as multivariate kernel densities. 
% https://en.wikipedia.org/wiki/Kernel_density_estimation
The bandwidth matrix $H$ has a huge impact on the accuracy of the estimator, as it is the most important hyperparameter for density smoothing. In the multidimensional scenario, both the width and the orientation are adjustable. The bandwidth has to be adapted to the considered state space. However, there are optimal  as well as rule-of-thumb criteria for the choice of $H$, see, e.\,g., \cite{Wand:1995}. A sketch of the proposed kernel density estimator of $\tilde{p}(x_{k-1}|Y_{k-1})$ in a multi-threaded environment is given in Algorithm 2. Note, however, that there are a number of machine learning packages that contain very efficient software implementations of the kernel density estimator, such as mlpack (C++ / Python) and sklearn (Python).

\begin{algorithm}[h]
	\caption{Kernel Density Estimation}
	\begin{algorithmic}[1]
		
		\State Set number of threads $N_{threads}$
		\State Split $\mathcal{A}$ in  $N_{threads}$ subsets $\mathcal{A}_1,\ldots,\mathcal{A}_{N_{threads}}$
		
		\ForAll{$\mathcal{A}_i \;\mathrm{in}\; \mathcal{A}_1,\ldots,\mathcal{A}_{N_{threads}}$}
		\State Start separate thread:
		\ForAll{$\alpha \;\mathrm{in}\; \mathcal{A}_i$}
		\State Evaluate coefficients (\ref{equ:approximated_coefficients}):
		\State{$m^\alpha_{k|k-1}  = \frac{ \sum_{i = 1}^{N_{s}} K\left( {H}^{-1/2} (x^{\alpha}_{k-1} - x_{k-1}^{(i)})\right)}{n \sqrt{\det({H})}}$}
		\EndFor
		\EndFor
		
	\end{algorithmic}
\end{algorithm}

%\subsection{Equations}

%Some words might be appropriate describing equation~(\ref{eq:sample}), if we had but time and space enough. 

%\begin{equation} \label{eq:sample}
%{{\partial F}\over {\partial t}} = D{{\partial^2 F}\over {\partial x^2}}.
%\end{equation}

%See \cite{Abl:56}, \cite{AbTaRu:54}, \cite{Keo:58} and \cite{Pow:85}.

%\subsubsection{Example.} This equation goes far beyond the
%celebrated theorem ascribed to the great Pythagoras by his followers.

%\begin{thm}   % use the thm environment for theorems
%The square of the length of the hypotenuse of a right triangle equals the sum of the squares of the lengths of the other two sides.
%\end{thm}

%\begin{pf}    % and the pf environment for proofs
%The square of the length of the hypotenuse of a right triangle equals the sum of the squares 
%of the lengths of the other two sides.
%\end{pf}

%% There are a number of predefined theorem-like environments in
%% ifacconf.cls:
%%
%% \begin{thm} ... \end{thm}            % Theorem
%% \begin{lem} ... \end{lem}            % Lemma
%% \begin{claim} ... \end{claim}        % Claim
%% \begin{conj} ... \end{conj}          % Conjecture
%% \begin{cor} ... \end{cor}            % Corollary
%% \begin{fact} ... \end{fact}          % Fact
%% \begin{hypo} ... \end{hypo}          % Hypothesis
%% \begin{prop} ... \end{prop}          % Proposition
%% \begin{crit} ... \end{crit}          % Criterion

%Of course LaTeX manages equations through built-in macros. You may wish to use the \texttt{amstex} package for enhanced math capabilities.

\subsection{Approximation of  the Inference Stage}

At the stage of inference at time $k$, the a priori density $\tilde{p}(x_k|Y_{k-1})$ is given as the current best estimate and a newly obtained measurement $y_k$ is available. These data are used in order to determine the a posteriori probability density $\tilde{p}(x_k|Y_k)$.
Since in the approximated state space, the pdfs are described by the coefficients $m^\alpha_{k|k}$  and $m^\alpha_{k|{k-1}}$, only their transition is considered.
Since the coefficients correspond more or less only to the evaluation of the densities at the grid nodes, the relationship is immediately given with (\ref{equ:bayesian_inference}) as
\begin{equation} \label{equ:approx_bayesian_inference}
m^\alpha_{k|k}  = \frac{p(y_k|x^\alpha_k)}{p(y_k|Y_{k-1})} m^\alpha_{k|k-1}
\end{equation} 
where $m^\alpha_{k|k} = p(x^\alpha_k|Y_{k})$ and $m^\alpha_{k|k-1} = p(x^\alpha_k|Y_{k-1})$. The denominator $p(y_k|Y_{k-1})$ corresponds to a normalization factor that does not affect the proportion of the distribution and the MAP estimation $\hat{x}^+_k$. Therefore, the normalization of the individual coefficients is omitted for the time being, since a subsequent scaling with the integral of the interpolated density ensures the normalization property. Since there are no dependencies between the coefficients $m^\alpha_{k|k}$ and $m^\beta_{k|k}$ for all $\alpha, \beta \in \mathcal{A}$ and $\alpha \neq \beta$, the set $\mathcal{A}$ can be split into several subset, which allows the coefficients to be calculated in parallel. The divide-and-conquer strategy aims to distribute the workload across multiple threads on a multi-core system to enable fast parallel execution. The numerator can be computed by evaluating (\ref{equ:bayesian_inferance_measurement_liklihood}) at the specific grid nodes $x_k^\alpha$. The calculation is not very complex from a computational point of view, and note that in the typical scenario of additive noise, even the Jacobian matrix becomes the identity matrix. A sketch of the proposed procedure is given in Algorithm 3.
\begin{algorithm}\label{algo:inference}
	\caption{Inference}
	\begin{algorithmic}[1]
		\setstretch{1.0}
		\State Set number of threads $N_{threads}$
		\State Split $\mathcal{A}$ in  $N_{threads}$ subsets $\mathcal{A}_1,\ldots,\mathcal{A}_{N_{threads}}$
		
		\ForAll{$\mathcal{A}_i \;\mathrm{in}\; \mathcal{A}_1,\ldots,\mathcal{A}_{N_{threads}}$}
		\State Start separate thread:
		\ForAll{$\alpha \;\mathrm{in}\; \mathcal{A}_i$}
		\State Evaluate numerator of (\ref{equ:bayesian_inferance_measurement_liklihood}):
		\State $\psi^\alpha_k = \left(p(v_k) \det(\partial_{y_k} g_k)^{-1}\right) \circ g_k(x^\alpha_k,u_k,y_k)$
		\State Evaluate coefficients (\ref{equ:approximated_coefficients}):
		\State{$m^\alpha_{k|k}  = \psi^\alpha_k\;m^\alpha_{k|k-1}$}
		\EndFor

		\EndFor
		
		\State Approximate density integral, see (\ref{equ:integral_approximation}):
		\State $ V_k = \sum_{i=1}^{N_{c}} \tilde{p}(x^i_c|Y_k) \Delta x^i_k$ 
		\State Normalize density:
		\ForAll{$\alpha \;\mathrm{in}\; \mathcal{A}$}
		\State{$m^\alpha_{k|k}  = V_k^{-1} m^\alpha_{k|k}$}
		\EndFor

	\end{algorithmic}
\end{algorithm}

\section{Experimental Results} \label{sec:experimental_results}
\begin{figure}[h]
	\begin{center}
		\includegraphics{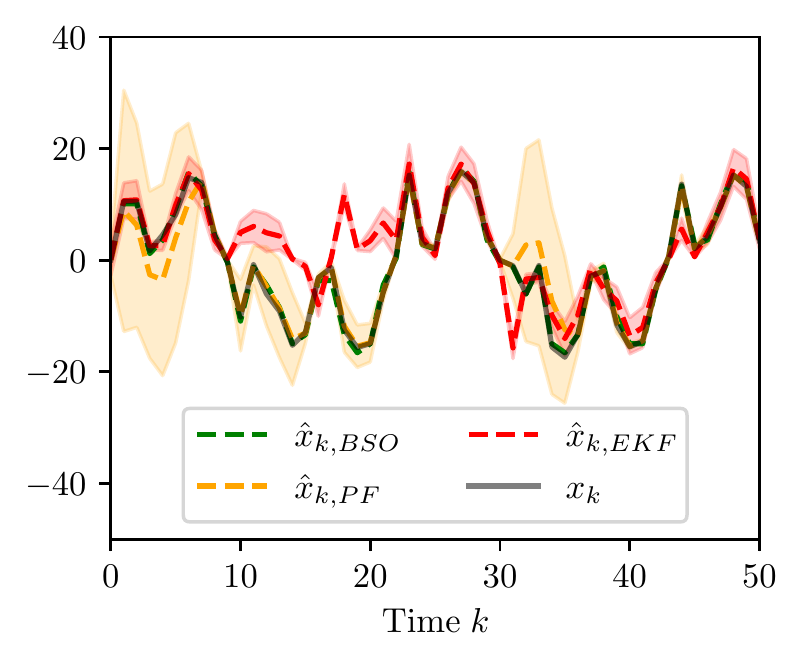}    % The printed column width is 8.4 cm.
		\caption{Benchmark Example: Comparison of the estimation performance of the extended Kalman filter (EKF), the particle filter (PF), and the Bayesian state observer (BSO) for the nonlinear system (\ref{equ:benchmark_example}).}
		\label{fig:benchmark_example_estimation}
	\end{center}
\end{figure}

\begin{figure}[ht]
	\begin{center}
		\includegraphics{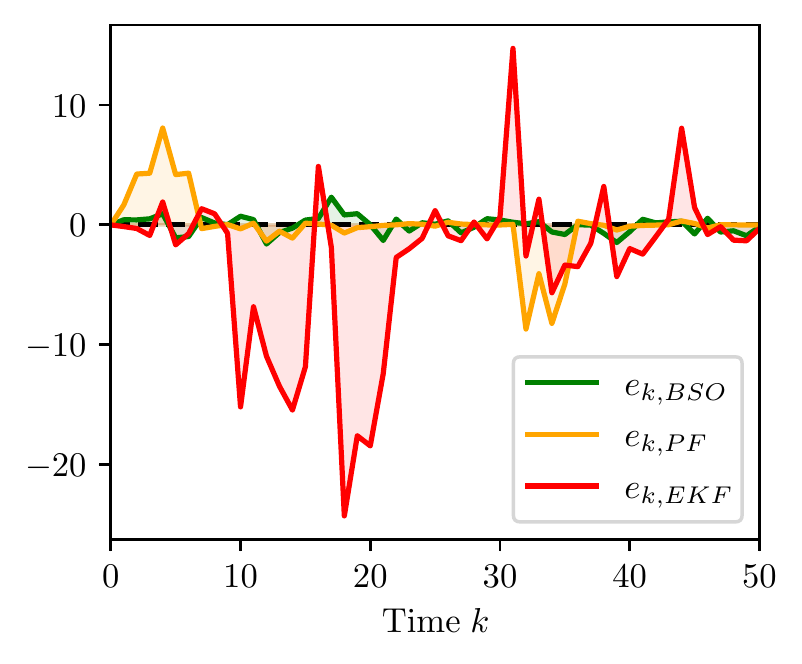}    % The printed column width is 8.4 cm.
		\caption{Benchmark Example: Comparison of the estimation error of the extended Kalman filter (EKF), the particle filter (PF), and the Bayesian state observer (BSO) for the nonlinear system (\ref{equ:benchmark_example}).} 
		\label{fig:benchmark_example_error}
	\end{center}
\end{figure}

%\subsection{Benchmark Example}
The previously presented algorithm is demonstrated on an autonomous scalar time-varying system with nonlinear state and measurement equation. 
The proposed test scenario is a benchmark in the nonlinear estimation literature, since the nonlinearity together with the time-dependent nature of the system turned out to be hardly capable by standard state observer.
The considered system with state $x_k \in \mathbb{R}$ and measured output $y_k \in \mathbb{R}$ is given by
\begin{subequations}
	\label{equ:benchmark_example}
	\begin{align}
	x_{k+1} &= \frac{1}{2} x_k + \frac{25 x_k}{1+ x^2_k} + 8 \cos(1.2 k) + w_k \\
	y_k &= \frac{1}{20} x^2_k + v_k
	\end{align}
\end{subequations}
where the additive process noise $w_k$ and measurement noise $v_k$ are modeled as zero-mean normal distributed random variables with variances equal to one, i.\,e., $w_k, v_k \sim \mathcal{N}(0,1)$.  The EKF and the PF are adopted from \cite{Sim:06}, whereas the Bayesian state observer corresponds to the previously presented grid-based implementation. Fig. \ref{fig:benchmark_example_estimation} demonstrates the estimation performance of the three proposed observer, where the PF is initialized with 500 particles and the BSO is implemented on a set of 500 grid points uniformly spaced between the boundaries -40 and 40. 
As Fig. \ref{fig:benchmark_example_estimation} shows, the EKF deviates far from the correct state in parts of the observation interval. What is remarkable, however, is that the correct state is sometimes even outside the observer's 95 percent confidence interval (red area). Although the confidence interval (orange area) of the PF outweighs that of the EKF in some places, the PF shows a significantly better performance. The estimates were only exceeded by the BSO, which is also clearly reflected in %Fig. \ref{fig:benchmark_example_error}
the following figure. Note, that the estimates depend on randomly generated samples. Therefore, between two series of experiments with different random seeds, the result may vary slightly. The linear estimation error $e_k = x_k - \hat{x}_k$ for the three proposed observer is illustrated in Fig. \ref{fig:benchmark_example_error}. The observer performance is assessed by the mean squared error $\text{MSE} = \frac{1}{n} \sum_{k=1}^{n} (x_k - \hat{x}_k)^2$ and CPU computing time for a estimation step. An overview of the performance with different numbers of cells, i.\,e. particles and grid points, is given in Tab. \ref{tb:observer_performance}, where the minimum and average performance values of one hundred simulations were recorded. The BSO outperformed the PF in both computation time and observer performance above a certain number of cells. Note that the PF is the standard implementation as proposed in \cite{Sim:06} and there are also extensions to improve the estimation performance. However, the hyperparameter of the BSO, i.\,e., grid width, grid limit as well as the number of grid points, were not optimized either. A significant increase in performance could still be achieved by adaptive regridding.

\begin{table}[h]
	\begin{center}
		\caption{Observer performance evaluated on a standard workstation (6 x Intel Core i5-85000 CPU @ 3.00 GHz and 8 GB RAM)}\label{tb:observer_performance}
		\begin{tabular}{lccrrcrr}
			\toprule
			
			\multirow{2}{*}{\vspace*{-0.3cm}Filter} &  \multirow{2}{*}{\vspace*{-0.3cm}Cells} &&  \multicolumn{2}{c}{Time (ms)} && \multicolumn{2}{c}{MSE} \\ 
			\cmidrule{4-5}
			\cmidrule{7-8}
			&  &&  avg.  & min.  & & avg.  & min. \\
			\midrule 
			 \vspace{0.1cm} 
			EKF & - && 0.04 & 0.03 &&  51.93& 20.77 \\
			PF & 100 && 11.16 & 10.09 && 25.16&8.16\\
			   & 200 && 22.64 & 21.24 && 19.32&7.10 \\
			 \vspace{0.1cm}
			& 500 && 75.46 & 73.46 && 14.92&5.78 \\
			
			BSO & 100 && 5.53 & 5.23 && 30.35 &10.35 \\ 
			& 200 && 8.83 & 7.44       && 13.02 &6.02 \\
			& 500 && 15.48 & 13.83     && 7.99  &0.45 \\
			\bottomrule
		\end{tabular}
	\end{center}
\end{table}

\section{Conclusion}

In this contribution, a piecewise linear approximation of the Bayesian state estimation problem for nonlinear systems has been presented. The approach proposes a Markov chain Monte Carlo propagation stage, where the Metropolis-Hastings algorithms is used for the sampling from the a posteriori density. The a priori density is reconstructed by means of a kernel density estimation of the propagated samples. The performance of the Bayesian state observer is verified on a benchmark example. The decisive advantage of this variant is its applicability for systems of higher order, which will be presented in a further paper.

%\begin{ack}
%Place acknowledgments here. 
%\end{ack}

\bibliography{ifacconf}             % bib file to produce the bibliography
                                                     % with bibtex (preferred)
                                                   
%\begin{thebibliography}{xx}  % you can also add the bibliography by hand

%\bibitem[Able(1956)]{Abl:56}
%B.C. Able.
%\newblock Nucleic acid content of microscope.
%\newblock \emph{Nature}, 135:\penalty0 7--9, 1956.

%\bibitem[Able et~al.(1954)Able, Tagg, and Rush]{AbTaRu:54}
%B.C. Able, R.A. Tagg, and M.~Rush.
%\newblock Enzyme-catalyzed cellular transanimations.
%\newblock In A.F. Round, editor, \emph{Advances in Enzymology}, volume~2, pages
%  125--247. Academic Press, New York, 3rd edition, 1954.

%\bibitem[Keohane(1958)]{Keo:58}
%\newblock \emph{Power and Interdependence: World Politics in Transitions}.
%\newblock Little, Brown \& Co., Boston, 1958.

%\bibitem[Powers(1985)]{Pow:85}
%\newblock Is there a way out?
%\newblock \emph{Harpers}, pages 35--47, June 1985.

%\bibitem[Soukhanov(1992)]{Heritage:92}
%\newblock \emph{{The American Heritage. Dictionary of the American Language}}.
%\newblock Houghton Mifflin Company, 1992.

%\end{thebibliography}

%\appendix
%\section{A summary of Latin grammar}    % Each appendix must have a short title.
%\section{Some Latin vocabulary}              % Sections and subsections are supported  
                                                                         % in the appendices.
                                                                                                                          
\end{document}